\documentclass[12pt]{article}%
\usepackage{graphicx}
\usepackage{amsmath}
\usepackage{amsfonts}
\usepackage{amssymb}%
\newtheorem{theorem}{Theorem}
\newtheorem{acknowledgement}[theorem]{Acknowledgement}

\newtheorem{conjecture}[theorem]{Conjecture}
\newtheorem{corollary}[theorem]{Corollary}

\newtheorem{lemma}[theorem]{Lemma}

\newtheorem{proposition}[theorem]{Proposition}
\newtheorem{remark}[theorem]{Remark}

\newenvironment{proof}[1][Proof]{\textbf{#1.} }{\ \rule{0.5em}{0.5em}}
\begin{document}

\title{TRANSVERSE RIEMANN-LORENTZ TYPE-CHANGING METRICS WITH POLAR END. }
\date{}
\author{J. Lafuente-L\'{o}pez\thanks{Departamento de Geometr\'{\i}a y Topolog\'{\i}a,
Facultad de Matem\'{a}ticas, Universidad Complutense, 28040 Madrid, Spain.
Email address: javier{\_}lafuente@mat.ucm.es}}
\maketitle

\begin{abstract}
Consider a smooth manifold $M$ with a smooth cometric $g^{\ast}$ which changes
the bilineal type by transverse way, on a hypersurface $D^{\infty}$. Suppose
that the radical annihilator hyperplane is tangent to $D^{\infty}$. We examine
the geometry of the ($g^{\ast}$-dual) covariant metric $g$ on $M-$ $D^{\infty
}$, prove the existence of a canonical (polar-normal) vectorfield whose
integral curves are $C^{\infty}-$pregeodesics crossing $D^{\infty}$
transversely for each point, and analyze the curvature behavior using a
natural coordinates. Finally we give an approach to the conformal geometry of
such spaces and suggest some application as cosmological big-bang model.


\end{abstract}

\section{Introduction}

There are several geometrical and physical reasons to study the metrics with
signature type changing Lorentz to Riemann (see for example the introductions
to \cite{kossow2}, \cite{kossow5} and \cite{kossow7}). For physical reasons,
two proposals for such spacetimes have been advanced:

a) The metric everywhere is smooth but it is degenerate at the hypersurface of
signature type change.

b) The metric is everywhere non degenerate but fails to be defined or to be
continuous at the hypersurface that divides the Riemannian from Lorentzian region.

There are many articles devoted to the proposal a) (see .\cite{kossow1}
,\cite{kossow2}, \cite{Hayward}, \cite{AgLa}) and some others to proposal b)
(see for example \cite{Ellis}, \cite{kossow7})

In this article we analyze a particular version\footnote{This version has been
inspired by a personal communication of Prof. M.Kossowski} of proposal b).
Here the dual metric $g^{\ast}$ (but not the metric $g$) is smooth and well
defined on the whole space $M$ and it is of transverse type changing from
Lorentz to Riemann trough a hypersurface $D^{\infty}$ (called polar). We refer
these metric $g$ (with certain annihilator condition added) as a
\emph{Lorentz-Riemann metric with polar end}. Its formal definition is
displayed at the beginning of the Section 2. The main result of this Section
is that the \emph{geometry }of these space allow to define canonically a
polar-normal transversal direction along $D^{\infty}$.

In Section 3 we prove the existence of an unique pregeodesic crossing
transversely $D^{\infty}$ for any $p\in D^{\infty}$. Moreover these
pregeodesics cross at the polar-normal direction. Then using as parameter of
the pregeodesics, the square of the arc length to $D^{\infty}$, we may
establish by a standard process, a natural coordinate $\left(  z_{1}%
,\ldots,z_{m}\right)  $ system around any point to $D^{\infty}$. The partial
$\partial_{z_{m}}$ is then a canonical polar-normal pregeodesic vectorfield.
Its flow are called \emph{the polar normal flow}$.$

Using the natural coordinates in Section 4, we analyze the behavior near of
$D^{\infty}$ of the semiriemannian curvatures. The conclusions are resumed in
Theorem \ref{Th2}.

The \emph{Lorentz-Riemann metric with polar end }character is preserved by a
conformal change. Section 5 is devoted to analyze some aspects of this
conformal geometry $\left(  M,\mathcal{C}\right)  $. The main result is that
fixing any (local) flow which moves $D^{\infty}$, there exist a metric
$g\in\mathcal{C}$ such that this flow is \emph{the polar normal flow }with
respect to $g.$ Moreover $g$ is univocally determined around $D^{\infty}$ .
Finally we consider the Lorentzian piece of $\left(  M,\mathcal{C}\right)  $
as the causal structure support of a admissible cosmological model where
$D^{\infty}$means the big-bang singularity, then we speculate with the
existence of a metric $g\in\mathcal{C}$, such that his \emph{ polar normal
flow }moves $D^{\infty}$ by (simultaneity) hypersurfaces with constant
sectional curvature. \newpage

\section{Preliminaries.}

\subsection{Type-changing metrics with polar end. \label{Sec1}}

A \emph{type-changing metric space with polar end } is a $m-$dimensional
manifold $M$ ( $m\geq2)$, endowed with a smooth, symmetric $(0,2)-$tensorfield
$g$ over an open set $M-D^{\infty}$of $M$. ($D^{\infty}\neq\varnothing$). For
any $p\in M-D^{\infty}$, we construct the dual metric at $p$, \emph{\ }%
$g^{\ast}:T_{p}^{\ast}M\times T_{p}^{\ast}M\rightarrow\mathbb{R}$, by
\begin{equation}
g^{\ast}\left(  \alpha,\beta\right)  =g\left(  X_{\alpha},X_{\beta}\right)
\label{g*}%
\end{equation}
and $g^{\ast}$ is a $\left(  2,0\right)  $ tensor over $M-D^{\infty}$. We demand:

\begin{itemize}
\item[D1)] \emph{The dual metric }$g^{\ast}$ of $g$ on $M-D^{\infty}$\emph{,
has smoothly extension to }$M$\emph{, and it is transverse type-changing over
}$D^{\infty}$.\newline The \emph{transverse type-changing property for the
}extension means that if $\left(  \theta^{1},\ldots,\theta^{m}\right)  $ is a
coframe on a neighborhood $U$ of $p\in D^{\infty}$, then $U\cap D^{\infty
}=\left\{  x\in U:\left.  \det\left(  g\left(  \theta^{a},\theta^{b}\right)
\right)  \right\vert _{x}=0\right\}  $, and $\det\left(  g^{\ast}\left(
\theta^{a},\theta^{b}\right)  \right)  =0$ is a local equation for $D^{\infty
}$ that is
\[
d_{x}\det\left(  g\left(  \theta^{a},\theta^{b}\right)  \right)  \neq0\text{,
}\forall x\in U\cap D^{\infty}\
\]
Of course this condition is frame-independent. In particular $D^{\infty}$ is
an hypersurface (called \emph{polar end}$)$ and at each point $p\in D^{\infty
}$ (called \emph{polar point}) the radical $Rad_{p}\left(  g^{\ast}\right)
\subset T_{p}^{\ast}M$ is one dimensional. Therefore the annihilator
\[
An\left(  Rad_{p}\left(  g^{\ast}\right)  \right)  =\left\{  X\in T_{p}%
M:\mu\left(  X\right)  =0\text{, }\forall\mu\in Rad_{p}\left(  g^{\ast
}\right)  \right\}
\]
is $\left(  m-1\right)  $-dimensional. Also the signature of $g^{\ast}$ (or
$g$) changes by $+1$ or $-1$ across $D^{\infty}$.(see \cite{kossow1} for details)
\end{itemize}

The last condition required is:

\begin{itemize}
\item[\textbf{D2)}] \emph{The annihilator is tangent to }$D^{\infty}%
$.\emph{\newline}That is, $An\left(  Rad_{p}\left(  g^{\ast}\right)  \right)
=T_{p}D^{\infty}$ for any $p\in D^{\infty}$. In fact if $\mu\in Rad_{p}\left(
g^{\ast}\right)  $, and $\mu\neq0$, the\emph{ }condition at $p$ say : $\ker
\mu=T_{p}D^{\infty}$.\medskip
\end{itemize}

We refer to $\left(  M,g\right)  $ as a type-changing metric space with
\emph{polar end} ($D^{\infty}$)

Of course we replace the term \emph{type-changing metric} by
\emph{Riemann-Lorentz}, if every component of $M-D^{\infty}$ is either Riemann
or Lorentz.

Henceforth we restrict our attention to the Riemann-Lorentz type.\medskip

If $N$ is manifold (possibly with boundary) denote $\mathfrak{X}\left(
N\right)  $ the $C^{\infty}(N)$-module of all vectorfields on $N$. If $D$ is
submanifold of $N,$ then $\mathfrak{X}_{N}\left(  D\right)  =\left\{
X\in\mathfrak{X}\left(  N\right)  :\left.  X\right\vert _{D}\in\mathfrak{X}%
\left(  D\right)  \right\}  $ is the submodule the vectorfields tangent to $D$.

Also $\mathfrak{X}_{D}\left(  N\right)  $ or ($\mathfrak{X}_{D}$ if we
understood $N$) is the $C^{\infty}(D)$ module of smooth $A:D\rightarrow TN$
with $A\left(  x\right)  \in T_{x}N$ for all $x\in D$. Finally $\Omega
^{1}\left(  N\right)  $ is the module of 1-forms (dual module of
$\mathfrak{X}\left(  N\right)  $)

We will use the following index conventions:$a,b,c\in\left\{  1,\ldots
,m\right\}  $ varies 1 to $m$, and $i,j,k\in\left\{  1,\ldots,m-1\right\}  $.
We also use Einstein's summation convention , unless the repeated index is $m$
We will work on a fixed neighborhood of a polar point $p\in D^{\infty}$ of the
Riemann-Lorentz space $\left(  M,g\right)  $ with polar end. Without loss
generality we will suppose that this neighborhood is the whole space $M.$ Also
we may suppose that $M-D^{\infty}$ have two connected component $D^{+}%
$(Riemannian) and $D^{-}$ (Lorentzian).

Let us consider some function $\tau\in C^{\infty}(M)$ such that $\,\tau
\mid_{D^{\infty}}\;=0$ and $\,d\tau\mid_{D^{\infty}}\;\neq0$ everywhere. We
say that \emph{$\tau=0$ is an equation for }$D^{\infty}$. Given another
function $f\in C^{\infty}(M),$ it holds: $\;f\mid_{D^{\infty}}%
\;=0\;\Leftrightarrow\;f=h\tau\,,$ for some $h\in C^{\infty}(M)$. When
$f\mid_{D^{\infty}}\;=0,$ we write $\;{\tau}^{-1}f\cong0\,$ and we say that
\emph{$\;\tau^{-1}f$ is extendible} as an element of $C^{\infty}(M)$.

\subsection{Polar-adapted frames.\label{Sec2}}

We say that a frame $\left(  E_{a}\right)  =\left(  E_{i},E_{m}\right)  $ on
$M$ is polar-adapted if $\left(  E_{i}\mid_{D^{\infty}}\right)  $ is a frame
of $D^{\infty}$, and%
\begin{equation}
\left(  g\left(  E_{a},E_{b}\right)  \right)  =\left(
\begin{array}
[c]{cc}%
\begin{tabular}
[c]{ccc}%
$1$ & $\ldots$ & $0$\\
$\vdots$ & $\ddots$ & $\vdots$\\
$0$ & $\ldots$ & $1$%
\end{tabular}
&
\begin{tabular}
[c]{c}%
$0$\\
$\vdots$\\
$0$%
\end{tabular}
\\%
\begin{tabular}
[c]{ccc}%
$0$ & $\ldots$ & $0$%
\end{tabular}
& 1/\tau
\end{array}
\right)  \label{gab}%
\end{equation}
where $\tau=0$ is a equation for $D^{\infty}$.

We say that the coframe $\left(  \theta^{1},\ldots,\theta^{m}\right)  $ on $M$
is $Rad^{\ast}-$adapted if%
\begin{equation}
\left(  g^{\ast}\left(  \theta^{a},\theta^{b}\right)  \right)  =\left(
\begin{array}
[c]{cc}%
\begin{tabular}
[c]{ccc}%
$1$ & $\ldots$ & $0$\\
$\vdots$ & $\ddots$ & $\vdots$\\
$0$ & $\ldots$ & $1$%
\end{tabular}
&
\begin{tabular}
[c]{c}%
$0$\\
$\vdots$\\
$0$%
\end{tabular}
\\%
\begin{tabular}
[c]{ccc}%
$0$ & $\ldots$ & $0$%
\end{tabular}
& \tau
\end{array}
\right)  \label{g*ab}%
\end{equation}
where $\tau=0$ is a equation for $D^{\infty}.$

If $\left(  \theta^{a}\right)  $ is coframe $Rad^{\ast}$-adapted and $\left(
E_{a}\right)  $ is the dual frame then $\left(  E_{a}\right)  $ is
polar-adapted frame. In fact (using notation of (\ref{alfa}) and (\ref{X})) we
have:
\begin{equation}
\delta_{i}^{a}=\theta^{a}\left(  E_{i}\right)  =g^{\ast}\left(  \theta
^{i},\theta^{a}\right)  =\theta^{a}\left(  X_{\theta^{i}}\right)  \label{dual}%
\end{equation}
since $\beta\left(  X_{\alpha}\right)  =g^{\ast}\left(  \alpha,\beta\right)
$. Analogously
\[
\left\{
\begin{tabular}
[c]{c}%
$\theta^{i}\left(  \tau E_{m}\right)  =0=g^{\ast}\left(  \theta^{i},\theta
^{m}\right)  =\theta^{i}\left(  X_{\theta^{m}}\right)  $\\
$\theta^{m}\left(  \tau E_{m}\right)  =\tau=g^{\ast}\left(  \theta^{m}%
,\theta^{m}\right)  =\theta^{m}\left(  X_{\theta^{m}}\right)  $%
\end{tabular}
\ \ \right.
\]

and we conclude that:%
\begin{equation}%
\begin{tabular}
[c]{cc}%
$X_{\theta^{i}}=E_{i}$ & $X_{\theta^{m}}=\tau E_{m}$\\
$\alpha_{E_{i}}=\theta^{i}$ & $\alpha_{E_{m}}=\left(  1/\tau\right)
\theta^{m}$%
\end{tabular}
\ \label{X_theta}%
\end{equation}

Taking account that $g\left(  X,Y\right)  =g^{\ast}\left(  \alpha_{X}%
,\alpha_{Y}\right)  $, in $M-D^{\infty}$ we get that
\[%
\begin{tabular}
[c]{c}%
$g\left(  E_{i},E_{j}\right)  =g^{\ast}\left(  \theta^{i},\theta^{j}\right)
=\delta_{j}^{i}$, $g\left(  E_{i},E_{m}\right)  =\left(  1/\tau\right)
g^{\ast}\left(  \theta^{i},\theta^{m}\right)  =0$\\
$g\left(  E_{m},E_{m}\right)  =\left(  1/\tau\right)  ^{2}g^{\ast}\left(
\theta^{m},\theta^{m}\right)  =1/\tau$%
\end{tabular}
\ \
\]
and .$\left(  g\left(  E_{a},E_{b}\right)  \right)  $ is as (\ref{gab})

Also $\theta^{m}\left(  E_{i}\right)  =0$, $\theta^{m}(p)\in Rad_{p}\left(
g^{\ast}\right)  $. Using the tangent annihilator property we conclude that
$\left(  \left.  E_{i}\right\vert _{D^{\infty}}\right)  $ is a frame for
$D^{\infty}$. Thus $\left(  E_{a}\right)  $ is a polar-adapted frame.

\begin{remark}
\label{Rem2}Since the existence of $Rad^{\ast}-$adapted coframes it is
straightforward, we conclude the existence of polar-adapted frames. On the
other and using a polar-adapted frame $\left(  E_{a}\right)  $ is easy to see
that $g$ induces by restriction on $D^{\infty}$ a canonical Riemannian
structure such that $\left(  \left.  E_{i}\right\vert _{D^{\infty}}\right)  $
is orthonormal frame.
\end{remark}

The same arguments show that if $\left(  E_{a}\right)  $ is a polar-adapted
frame then the dual $\left(  \theta^{a}\right)  $ is a $Rad^{\ast}$-adapted
coframe. We will prove now that it is possible to make a $Rad^{\ast}-$adapted
coframe $\left(  \theta^{i},\theta^{m}\right)  $ for any fixed $\mu=\theta
^{m}$:

\begin{proposition}
\label{mu}Let $\mu$ be a $1-$form on $M$ which is $Rad^{\ast}-$adapted (that
is $\mu\left(  x\right)  \in Rad_{x}\left(  g^{\ast}\right)  -\left\{
0\right\}  $ for all $x\in D^{\infty.}$). Then there exist (locally) a
$Rad^{\ast}$ $-$adapted coframe $\left(  \mu^{i},\mu^{m}=\mu\right)  $.
\end{proposition}

\begin{proof}
We start with an auxiliary $Rad^{\ast}$- adapted coframe $\left(  \theta
^{a}\right)  $ as in (\ref{g*ab}) Without lost generality we may write
\[
\mu=\sum\left(  \tau h_{i}\right)  \theta^{i}+\theta^{m}%
\]
since if $\left(  \mu^{i},\mu^{m}=\mu\right)  $ is $Rad^{\ast}$- adapted then
the same holds for $\left(  \mu^{i},h\mu\right)  $ for any smooth everywhere
non-null $h$. Thus $\alpha=\sum X_{a}\theta^{a}$ is orthogonal to $\mu$ iff
$\sum h_{j}X_{j}+X_{m}=0$, and the co-distribution $\mu^{\perp}$ is generated
by the $\left(  m-1\right)  $-coframe $\left(  \vartheta^{i}=\theta^{i}%
-h_{i}\theta^{m}\right)  $ which is well-defined also over $D^{\infty}$. By
application of the the classical orthonormalization Graham Smith process to
$\left(  \vartheta^{i}\right)  $ we obtain the desired coframe.
\end{proof}

The dual result is the following:

\begin{theorem}
\label{Th1}Let $N$ be any vectorfield on $M$ transversal to $D^{\infty}$. Then
there exist (locally) a polar-adapted frame $\left(  N_{i},N_{m}=N\right)  $.
\end{theorem}

\begin{proof}
We start with an auxiliary $Rad^{\ast}$- adapted coframe $\left(  \theta
^{a}\right)  $ with (polar-adapted) dual frame, $\left(  E_{a}\right)  $. Thus
$\left(  g_{ab}\right)  $ is as (\ref{gab}). Writing $N=\sum h^{a}E_{a}$,
transversality implies that $h^{m}\left(  x\right)  \neq0$ for all $x\in
D^{\infty}$.Now taking $\mu=\sum\left(  \tau h^{i}\right)  \theta^{i}%
+h^{m}\theta^{m}$ we have:%
\begin{align}
X_{\mu}  &  =\sum\tau h^{i}X_{\theta^{i}}+h^{m}X_{\theta^{m}}\label{Xmu}\\
&  =\sum\tau h^{i}E_{i}+\tau h^{m}E_{m}\nonumber\\
&  =\tau N\nonumber
\end{align}
\newline On the other hand we may construct $\left(  \mu^{i},\mu^{m}%
=\mu\right)  $ coframe $Rad^{\ast}$-adapted as in proposition \ref{mu}.Let
$\left(  N_{i},N_{m}\right)  $ be his polar-adapted dual coframe. Then by
(\ref{X_theta}) is $X_{\mu}=\tau_{\mu}N_{m}$ where $\tau_{\mu}=g^{\ast}\left(
\mu,\mu\right)  $. Since $\left(  \tau_{\mu}=0\right)  $ is equation of
$D^{\infty}$ then there exist a smooth everywhere non-null $h$ such that
$\tau_{\mu}=h\tau.$ Therefore $X_{\mu}=\tau hN_{m}$. Comparing with
(\ref{Xmu}) we get that $N=hN_{m}$. But if $\left(  N_{i},N_{m}\right)  $ is
polar-adapted, the same is true for $\left(  N_{i},h^{-1}N_{m}=N\right)  $.

\begin{corollary}
\label{Cor1}Given any vectorfield $N$ on $M$ transversal to $D^{\infty}$ and
$X\in\mathfrak{X}_{M}\left(  D^{\infty}\right)  $, we have $g\left(
X,N\right)  \cong0$
\end{corollary}
\end{proof}

\begin{proof}
By the theorem there exist a polar-adapted frame $\left(  E_{i},E_{m}%
=N\right)  $ as in (\ref{gab})$.$Therefore we may write $X=\sum X^{i}%
E_{i}+\tau hE_{m}$ where $h$ is some smooth function on $M$. Then $g\left(
X,N\right)  =\sum X^{i}g_{im}+h$ which it is a differentiable function.
\end{proof}

\begin{remark}
\label{Rem1}Using polar-adapted frames $\left(  E_{a}\right)  $ (as in
(\ref{gab})) is easily to prove that $g\left(  X,Y\right)  \cong0$ if $X$ or
$Y$ belongs to $\mathfrak{X}_{M}\left(  D^{\infty}\right)  $. In fact if
$X=\sum X^{a}E_{a}$ ...etc. then$\left.  X^{m}Y^{m}\right\vert _{D^{\infty}%
}=0$ thus $g\left(  X,Y\right)  =\sum X^{i}Y^{i}+\tau^{-1}\left(  X^{m}%
Y^{m}\right)  \in C^{\infty}\left(  M\right)  $. \newline On the other hand
note that $\tau g$ is defined on the whole space $M$, if $\left(
\tau=0\right)  $ is an equation for $D^{\infty}$.
\end{remark}

\subsection{The dual connection near to D$^{\infty}$.\label{Sec21}}

First we remark the local nature of the work. In fact we should be replaced in
any case $M$ by a suitable neighborhood of a polar point $p\in D^{\infty}$.

We recall (see \cite{kossow3} for more details) that on the Riemann-Lorentz
space $(M-D^{\infty},g)$ there exists a unique torsion-free metric \emph{dual
connection}, which it is characterized as the unique map $\square
:\mathfrak{X}(M-D^{\infty})\times\mathfrak{X}(M-D^{\infty})\rightarrow
\mathfrak{X}^{\ast}(M-D^{\infty})$ satisfying, for all $A,B,C\in
\mathfrak{X}(M-D^{\infty})$, the Koszul-like formula:%

\begin{equation}
\left.
\begin{array}
[c]{l}%
2\square_{A}B(C):=A\left\langle B,C\right\rangle +B\left\langle
C,A\right\rangle -C\left\langle A,B\right\rangle \\
\;\;\;\;\;\;\;\;\;\;\;-\left\langle A,[B,C]\right\rangle +\left\langle
B,[C,A]\right\rangle +\left\langle C,[A,B]\right\rangle \;\;.
\end{array}
\right.  \label{square}%
\end{equation}

It follows that $\square$ is compatible with the Levi-Civita connection
$\nabla$ on $M-\left(  D^{0}\cup D^{\infty}\right)  $, in the sense that it
holds: $\,\square_{A}B(C)=\left\langle \nabla_{A}B,C\right\rangle $.

With respect to the frame $\left(  E_{i},E_{m}\right)  $ the dual connection
is determined by the Christopher symbols $\Gamma_{cab}=\square_{E_{a}}%
E_{b}(E_{c}):$%
\begin{equation}
\left(  \square_{X}Y\right)  (Z)=X(Y^{b})g_{bc}Z^{c}+\Gamma_{cab}Z^{c}%
X^{a}Y^{b}\text{ .} \label{dualZXY}%
\end{equation}
for $g_{ab}=g(E_{a},E_{b)}$, $X=X^{a}E_{a}$,...etc. Explicitly:%
\begin{equation}
\Gamma_{cab}=\frac{1}{2}\left\{
\begin{tabular}
[c]{c}%
$E_{a}\left(  g_{bc}\right)  \quad+\quad E_{b}\left(  g_{ca}\right)
\quad-\quad E_{c}\left(  g_{ab}\right)  $\\
$-g\left(  E_{a},\left[  E_{b},E_{c}\right]  \right)  +g\left(  E_{b},\left[
E_{c},E_{a}\right]  \right)  +g\left(  E_{c},\left[  E_{a},E_{b}\right]
\right)  $%
\end{tabular}
\ \ \ \ \ \ \right\}  \label{Gamm}%
\end{equation}
if the frame is polar-adapted then $\left(  g_{ab}\right)  $ is as
(\ref{gab}). In particular note that
\begin{equation}
\Gamma_{kij}\cong0 \label{Gamm_kij}%
\end{equation}
since $\left[  E_{j},E_{k}\right]  \in\mathfrak{X}_{M}\left(  D^{\infty
}\right)  $ and (by Remark \ref{Rem1}) $g\left(  E_{i},\left[  E_{j}%
,E_{k}\right]  \right)  \cong0$. Moreover%
\[
\tau\Gamma_{mij}=\frac{1}{2}\left\{  -\tau g\left(  E_{i},\left[  E_{j}%
,E_{m}\right]  \right)  +\tau g\left(  E_{j},\left[  E_{m},E_{i}\right]
\right)  +\tau g\left(  E_{m},\left[  E_{i},E_{j}\right]  \right)  \right\}
\]
writing $\left[  E_{a},E_{b}\right]  =\sum C_{ab}^{c}E_{c}$, and taking
account the look of $\left(  \tau g_{ab}\right)  $ we obtain%
\[%
\begin{tabular}
[c]{c}%
$\tau g\left(  E_{i},\left[  E_{j},E_{m}\right]  \right)  =\tau g\left(
E_{i},\sum C_{jm}^{a}E_{a}\right)  =\tau C_{jm}^{i}$\\
$\tau g\left(  E_{m},\left[  E_{i},E_{j}\right]  \right)  =\tau g\left(
E_{m},\sum C_{ij}^{a}E_{a}\right)  =C_{ij}^{m}$%
\end{tabular}
\ \
\]
but $\left.  C_{ij}^{m}\right\vert _{D^{\infty}}=0$ (since $\left[
E_{i},E_{j}\right]  \in\mathfrak{X}_{M}\left(  D^{\infty}\right)  $. Therefore
$\left.  \tau\Gamma_{mij}\right\vert _{D^{\infty}}=0$ and
\begin{equation}
\Gamma_{mij}\cong0\,\ \text{analogous }\Gamma_{kmj}\cong0\text{, }\Gamma
_{kim}\cong0 \label{Gamm_mij}%
\end{equation}
because (again by Remark \ref{Rem1}) $\tau g$ is defined on the whole $M$.
Taking account
\[
\tau E_{k}\left(  \frac{1}{\tau}\right)  =-\frac{E_{k}\left(  \tau\right)
}{\tau}\cong0\text{ (since}\left.  E_{k}\left(  \tau\right)  \right\vert
_{D^{\infty}}\text{ }=0\text{)}%
\]
we conclude that
\begin{align}
&  \frac{1}{2}\left\{  -\tau E_{k}\left(  \frac{1}{\tau}\right)  -2\tau
g\left(  E_{m},\left[  E_{m},E_{k}\right]  \right)  \right\} \nonumber\\
&  =\tau\Gamma_{kmm}\cong0\text{ analogous }\tau\Gamma_{mim}\cong0\text{,
}\tau\Gamma_{mmj}\cong0 \label{Gamm_kmm}%
\end{align}

However $\left.  E_{m}\left(  \tau\right)  \right\vert _{D^{\infty}}$ is non
null everywhere and%
\begin{equation}
\frac{1}{2}\tau^{2}E_{m}\left(  \frac{1}{\tau}\right)  =-\frac{E_{m}\left(
\tau\right)  }{2}=\tau^{2}\Gamma_{mmm} \label{Gamm_mmm}%
\end{equation}

A first consequence of these computations are:

\begin{proposition}
Let $X,Y,Z\in\mathfrak{X}\left(  M\right)  $ then:

\begin{enumerate}
\item If two of these tree vectorfields are tangent to $D^{\infty}$ then
$\left(  \square_{Z}X\right)  (Y)\cong0$ .

\item If $Z$ is transversal to $D^{\infty}$, $X\in\mathfrak{X}_{M}\left(
D^{\infty}\right)  $ and $g\left(  X,Z\right)  =0$ then (locally)
$g(X,X)^{-1}\left(  \square_{Z}Z\right)  (X)\cong0$

\item If $Z$ is transversal to $D^{\infty}$ then for any $V\in\mathfrak{X}%
\left(  D^{\infty}\right)  $ the function%
\begin{equation}
\beta_{Z}\left(  V\right)  =\left.  \frac{\left(  \square_{Z}Z\right)
(X)}{g(Z,Z)}\right\vert _{D^{\infty}}\text{ }\left(  \left.  X\right\vert
_{D^{\infty}}=V\right)  \label{beta}%
\end{equation}
is independent of $X\in\mathfrak{X}_{M}\left(  D^{\infty}\right)  $ such that
$\left.  X\right\vert _{D^{\infty}}=V$ and $g\left(  X,Z\right)  =0$. Moreover
$\beta_{Z}\in\Omega^{1}\left(  D^{\infty}\right)  $
\end{enumerate}
\end{proposition}

\begin{proof}
First we take an auxiliary polar-adapted frame $\left(  E_{a}\right)  $ as in
(\ref{gab}).

The assert 1. is an easy consequence of the formula (\ref{dualZXY}), using the
general expression $X=\sum X^{i}E_{i}+\tau hE_{m}$, of any $X\in
\mathfrak{X}_{M}\left(  D^{\infty}\right)  $, and taking account that $\tau
g$, $\tau\Gamma_{kmm}$, $\tau\Gamma_{mij}$,...and $\tau^{2}\Gamma_{mmm}$ are
defined on the whole space $M$.

To prove 2.and 3 note that by the Theorem \ref{Th1} we may suppose as well
(without lost generality) that $Z=E_{m}$. Then $Z^{m}=1,$ $0=Z^{i}=X^{m}$ and
$g(Z,Z)^{-1}=\tau$. Applying (\ref{dualZXY}) gives
\begin{equation}
\frac{\left(  \square_{Z}Z\right)  (X)}{g(Z,Z)}=\sum\tau\Gamma_{kmm}X^{i}%
\cong0 \label{beta1}%
\end{equation}
by (\ref{Gamm_kmm}). On the other hand if $\left.  \tau\Gamma_{kmm}\right\vert
_{D^{\infty}}=\gamma_{k}$, we write the 1-form claimed in (\ref{beta}):%
\begin{equation}
\beta_{E_{m}}=\sum\gamma_{k}\theta_{\infty}^{k}\text{ where }\gamma
_{k}=\left.  \tau\Gamma_{kmm}\right\vert _{D^{\infty}} \label{beta2}%
\end{equation}
where $\left(  \theta_{\infty}^{k}\right)  $ is the dual coframe of $\left(
\left.  E_{k}\right\vert _{D^{\infty}}\right)  .$This proves 3.
\end{proof}

We will require the following

\begin{lemma}
\label{Lem1}For any generic transversal vectorfield $Z\in\mathfrak{X}\left(
M\right)  $ and any smooth ever non null function $h$ we have $\beta_{Z}%
=\beta_{hZ}$.
\end{lemma}

\begin{proof}
Let $X\in$ $\mathfrak{X}\left(  D^{\infty}\right)  $ be a vectorfield and let
$\overline{X}\in\mathfrak{X}\left(  M\right)  $ be such that $\left.
\overline{X}\right\vert _{D^{\infty}}=X$ and $g\left(  \overline{X},Z\right)
=0$ (thus $g\left(  \overline{X},hZ\right)  =0$). We have%
\begin{align*}
\beta_{hZ}\left(  X\right)   &  =\left.  \frac{\left(  \square_{hZ}hZ\right)
(\overline{X})}{g(hZ,hZ)}\right\vert _{D^{\infty}}\\
&  =\left.  \frac{hZ\left(  h\right)  g\left(  \overline{X},Z\right)
+h^{2}\left(  \square_{Z}Z\right)  \left(  \overline{X}\right)  }{h^{2}%
g(Z,Z)}\right\vert _{D^{\infty}}\\
&  =\left.  \frac{\left(  \square_{Z}Z\right)  (\overline{X})}{g(Z,Z)}%
\right\vert _{D^{\infty}}=\beta_{Z}\left(  X\right)
\end{align*}

\end{proof}

\subsection{Polar-normal vectorfield.\label{Sec22}}

We say that the vectorfield $Z$ on $M$ is a polar-normal vectorfield if it is
transversal to $D^{\infty}$ and the associated $1-$form $\beta_{Z}$ on
$D^{\infty}$ is identically null. In order to prove the existence of
polar-normal vectorfield, we start with an auxiliary polar-adapted frame
$\left(  E_{i},E_{m}\right)  $ as (\ref{gab}) and let $\left(  E_{i}^{\infty
}\right)  =$ $\left(  \left.  E_{i}\right\vert _{D^{\infty}}\right)  $ be the
restriction frame on $D^{\infty}$. As in (\ref{beta2}) we have
\begin{equation}
\beta_{E_{m}}\left(  E_{k}^{\infty}\right)  =\gamma_{k}=\left.  \tau
\Gamma_{kmm}\right\vert _{D^{\infty}} \label{gamma_k}%
\end{equation}

We find a transversal to $D^{\infty}$ vectorfield
\begin{equation}
\overline{E}_{m}=\sum\lambda^{i}E_{i}+E_{m}=\widetilde{E}_{m}+E_{m}%
\in\mathfrak{X}\left(  M\right)  \text{ (where }g\left(  \widetilde{E}%
_{m},E_{m}\right)  =0\text{)} \label{E_m}%
\end{equation}
such that $\beta_{\overline{E}_{m}}\left(  E_{k}^{\infty}\right)  =0$. By the
Theorem \ref{Th1}, There exist $\left(  \overline{E}_{i}\right)  $ such that
$\left(  \overline{E}_{i},\overline{E}_{m}\right)  $ is polar-adapted frame,
and we may suppose without lost generality that $\left.  \overline{E}%
_{i}\right\vert _{D^{\infty}}=E_{i}^{\infty}$. (If not we may find an
orthogonal functional matrix $\left(  a_{j}^{i}\right)  $ with $a_{j}^{i}\in
C^{\infty}\left(  D^{\infty}\right)  $ such that $E_{j}^{\infty}=\sum
a_{j}^{i}\left.  \overline{E}_{i}\right\vert _{D^{\infty}}$ and we make
$\left(  \widehat{E}_{i},\overline{E}_{m}\right)  $ a polar-adapted frame,
with $\widehat{E}_{j}=\sum A_{j}^{i}E_{i}$ where $A_{j}^{i}\in C^{\infty
}\left(  M\right)  $ are the unique smooth function such that $\left.
A_{j}^{i}\right\vert _{D^{\infty}}=a_{j}^{i}$ and $\overline{E}_{m}\left(
A_{j}^{i}\right)  =0$). Therefore we may write for some smooth $\mu_{j}$ and
$\widetilde{E}_{i}$
\begin{equation}
\overline{E}_{i}=\widetilde{E}_{i}+\tau\mu_{i}E_{m}\text{ (where }g\left(
\widetilde{E}_{i},E_{m}\right)  =0\text{, and }\left.  \widetilde{E}%
_{i}\right\vert _{D^{\infty}}=E_{i}^{\infty}\text{)} \label{E_i}%
\end{equation}

Taking account (\ref{E_m}), (\ref{E_i}) and (\ref{gab}) we get%
\[
0=g\left(  \overline{E}_{i},\overline{E}_{m}\right)  =g\left(  \widetilde
{E}_{i},\widetilde{E}_{m}\right)  +\mu_{i}%
\]

Since $\left.  \widetilde{E}_{m}\right\vert _{D^{\infty}}=\sum g\left(
\widetilde{E}_{m},E_{i}^{\infty}\right)  E_{i}^{\infty}=\sum g\left(
\widetilde{E}_{m},\widetilde{E}_{i}\right)  E_{i}^{\infty}=\sum-\left.
\mu_{i}E_{i}\right\vert _{D^{\infty}}$, therefore%
\begin{equation}
\left.  \lambda^{i}\right\vert _{D^{\infty}}=-\left.  \mu_{i}\right\vert
_{D^{\infty}} \label{landamu}%
\end{equation}

Now we compute $\beta_{\overline{E}_{m}}\left(  E_{k}^{\infty}\right)  $ in
order to find $\lambda^{i}$ 'on $D^{\infty}$ which makes these values
identically null.

Taking account (\ref{E_m}) we have: $\square_{\overline{E}_{m}}\overline
{E}_{m}=\square_{\widetilde{E}_{m}}\widetilde{E}_{m}+\square_{\widetilde
{E}_{m}}E_{m}+\square_{E_{m}}\widetilde{E}_{m}+\square_{E_{m}}E_{m}$, where

\begin{itemize}
\item $\square_{\widetilde{E}_{m}}\widetilde{E}_{m}\cong0$ by (\ref{Gamm_kij})

\item $\square_{\widetilde{E}_{m}}E_{m}=\sum\lambda^{i}\square_{E_{i}}%
E_{m}=\sum\lambda^{i}\Gamma_{kim}\theta^{k}+\sum\lambda^{i}\Gamma_{mim}%
\theta^{m}$ (where $\left(  \theta^{a}\right)  $ is the dual coframe of
$\left(  E_{a}\right)  $)

\item $\square_{E_{m}}\widetilde{E}_{m}=\square_{E_{m}}\left(  \sum\lambda
^{i}E_{i}\right)  =\sum E_{m}\left(  \lambda^{i}\right)  \theta^{i}%
+\sum\lambda^{i}\Gamma_{kmi}\theta^{k}+\sum\lambda^{i}\Gamma_{mmi}\theta^{m}$
\end{itemize}

and we conclude using (\ref{Gamm_mij}) that for some $\theta\in\Omega
^{1}\left(  M\right)  $:
\[
\square_{\overline{E}_{m}}\overline{E}_{m}=\theta+\sum\lambda^{i}\left(
\Gamma_{mim}+\Gamma_{mmi}\right)  \theta^{m}+\square_{E_{m}}E_{m}%
\]
and taking account (\ref{beta}) and that $\overline{E}_{k}$ is extension of
$E_{k}^{\infty}$ which it is $g$-orthogonal to $\overline{E}_{m}$ we have:
\begin{align*}
\beta_{\overline{E}_{m}}\left(  E_{k}^{\infty}\right)   &  =\left.
\frac{\square_{\overline{E}_{m}}\overline{E}_{m}}{g\left(  \overline{E}%
_{m},\overline{E}_{m}\right)  }\left(  \overline{E}_{k}\right)  \right\vert
_{D^{\infty}}\\
&  =\left.  \frac{\tau\square_{\overline{E}_{m}}\overline{E}_{m}}{1+\tau
\sum\left(  \lambda^{i}\right)  ^{2}}\left(  \widetilde{E}_{k}+\tau\mu
_{k}E_{m}\right)  \right\vert _{D^{\infty}}\\
&  =\left.  \frac{\tau\left(  \theta+\sum\lambda^{i}\left(  \Gamma
_{mim}+\Gamma_{mmi}\right)  \theta^{m}\right)  +\tau\square_{E_{m}}E_{m}%
}{1+\tau\sum\left(  \lambda^{i}\right)  ^{2}}\left(  \widetilde{E}_{k}+\tau
\mu_{k}E_{m}\right)  \right\vert _{D^{\infty}}%
\end{align*}
but

\begin{itemize}
\item $\theta^{m}\left(  \widetilde{E}_{k}\right)  =0$ (by duality),

\item $\tau^{2}\left.  \left(  \Gamma_{mim}+\Gamma_{mmi}\right)  \right\vert
_{D^{\infty}}=0$ and $\tau^{2}\left.  \square_{E_{m}}E_{m}\left(
\widetilde{E}_{k}\right)  \right\vert _{D^{\infty}}=0$ (by (\ref{Gamm_kmm})),

\item $\tau\left.  \square_{E_{m}}E_{m}\left(  \widetilde{E}_{k}\right)
\right\vert _{D^{\infty}}=\gamma_{k}$ (see (\ref{gamma_k})) and

\item $\tau^{2}\square_{E_{m}}E_{m}\left(  E_{m}\right)  =-E_{m}\left(
\tau\right)  /2$ (by (\ref{Gamm_mmm})), thus:
\end{itemize}

\begin{equation}
\beta_{\overline{E}_{m}}\left(  E_{k}^{\infty}\right)  =\gamma_{k}-\frac{1}%
{2}\left(  \left.  \mu_{i}\right\vert _{D^{\infty}}\right)  \left.
E_{m}\left(  \tau\right)  \right\vert _{D^{\infty}} \label{beta3}%
\end{equation}
and if we want $\beta_{\overline{E}_{m}}\left(  E_{k}^{\infty}\right)  =0$ we
must to select $\lambda^{i}$ such that (see(\ref{landamu}))
\begin{equation}
\left.  \lambda^{i}\right\vert _{D^{\infty}}=-\frac{2\gamma_{k}}{\left.
E_{m}\left(  \tau\right)  \right\vert _{D^{\infty}}}\text{ \ (where }%
\gamma_{k}=\left.  \tau\Gamma_{kmm}\right\vert _{D^{\infty}}\text{)}
\label{landa_i}%
\end{equation}

We are now ready to prove the following main theorem:

\begin{theorem}
There exist $N\in\mathfrak{X}\left(  M\right)  $ polar-normal vectorfield
(that is $N$ is transversal to $D^{\infty}$ and $\beta_{N}=0$). Moreover if
$\overline{N}$ is other polar-normal vectorfield, then $N$ and $\overline{N}$
are proportional along $D^{\infty}$. (Thus a \emph{polar-normal direction}
along $D^{\infty}$is canonically determined ).
\end{theorem}

\begin{proof}
We have already proved the existence of polar-normal vectorfield. Only the
second part need be to prove:

Continuing with the previous argument let $N=E_{m}$ be the first polar normal
(that is $\beta_{E_{m}}=0$, thus $\gamma_{k}=0)$. Using Lemma \ref{Lem1} we
may suppose without lost generality that the other polar-normal $\overline
{N}=\overline{E}_{m}=$ $\widetilde{E}_{m}+E_{m}$ is as in (\ref{E_m}). Then by
(\ref{beta3}) and (\ref{landamu}) we have:
\[
0=\beta_{\overline{E}_{m}}\left(  E_{k}^{\infty}\right)  =-\frac{1}{2}\left(
\left.  \lambda^{i}\right\vert _{D^{\infty}}\right)  \left.  E_{m}\left(
\tau\right)  \right\vert _{D^{\infty}}%
\]
and this implies that $\left.  \lambda^{i}\right\vert _{D^{\infty}}=0,$
$\widetilde{E}_{m}=0$ and $\left.  \overline{N}\right\vert _{D^{\infty}%
}=\left.  \overline{N}\right\vert _{D^{\infty}}$.
\end{proof}

\section{Polar-adapted coordinates.\label{Sec3}}

The objective of this section is to prove that there exist a (essentially
unique) $C^{\infty}$-pregeodesic line traversing $D^{\infty}$ for everyone of
their points. Moreover each one traverse in the polar normal direction. The
proof is analogous to the proof in \cite{kossow1} \ of the similar result, in
\ the singular context.

Next using these pregeodesics we construct an special $C^{\infty}%
$-polar-adapted coordinates neighboring each point of $D^{\infty}$.

\subsection{Polar-normal pregeodesic.\label{Sec31}}

We start with a polar-normal adapted frame. This means a polar adapted frame
$\left(  E_{a}\right)  $ as in (\ref{gab}) such that $E_{m}$ is polar-normal.
In particular if $\Gamma_{cab}=\square_{E_{a}}E_{b}(E_{c})$ are the
Christopher symbols, we make the \emph{other }Christopher symbols $\Gamma
_{ab}^{c}$ defined by $\nabla_{E_{a}}E_{b}=\Gamma_{ab}^{c}E_{c}$ as:
\[
\left(
\begin{array}
[c]{c}%
\begin{tabular}
[c]{c}%
$\Gamma_{ab}^{1}$\\
\\
$\Gamma_{ab}^{m-1}$%
\end{tabular}
\\
\Gamma_{ab}^{m}%
\end{array}
\right)  =\left(
\begin{array}
[c]{cc}%
\begin{tabular}
[c]{ccc}%
$1$ & $\ldots$ & $0$\\
$\vdots$ & $\ddots$ & $\vdots$\\
$0$ & $\ldots$ & $1$%
\end{tabular}
&
\begin{tabular}
[c]{c}%
$0$\\
$\vdots$\\
$0$%
\end{tabular}
\\%
\begin{tabular}
[c]{ccc}%
$0$ & $\ldots$ & $0$%
\end{tabular}
& \tau
\end{array}
\right)  \left(
\begin{array}
[c]{c}%
\begin{tabular}
[c]{c}%
$\Gamma_{1ab}$\\
\\
$\Gamma_{m-1,ab}$%
\end{tabular}
\\
\Gamma_{mab}%
\end{array}
\right)
\]

This symbols controls the Levi-Civita Connection by the formula%
\[
\nabla_{X}Y=\left\{  X^{a}E_{a}(Y^{c})+\Gamma_{ab}^{c}X^{a}Y^{b}\right\}
E_{c}%
\]
\bigskip

Since $E_{m}$ is polar-normal we have $\gamma_{k}=\left.  \tau\Gamma
_{kmm}\right\vert _{D^{\infty}}=0$ and%
\[
\left.  \tau\Gamma_{kmm}\right\vert _{D^{\infty}}=\left.  \tau\Gamma_{mm}%
^{k}\right\vert _{D^{\infty}}=0
\]

Recollecting the information of (\ref{Gamm_kij}), (\ref{Gamm_kmm}),
(\ref{Gamm_mij}) and (\ref{Gamm_mmm}) we obtain:%
\begin{equation}%
\begin{tabular}
[c]{l}%
$\Gamma_{ij}^{c}\cong0$, $\left.  \Gamma_{ij}^{m}\right\vert _{D^{\infty}}%
=0$\\
$\Gamma_{mj}^{c}\cong0,$ $\Gamma_{im}^{c}\cong0$\\
$\tau\Gamma_{mm}^{m}=-\frac{1}{2}E_{m}\left(  \tau\right)  $%
\end{tabular}
\ \ \ \label{Gammas}%
\end{equation}

Now we follow analogous argument that in Theorem 2 of \cite{kossow1}: \bigskip

We fix a coordinate system $\left(  x^{i},x^{m}\right)  $ which we suppose
global (without loss generality), and take on $TM$ mixed coordinates $\left(
x^{a},u^{a}\right)  $, such that the following hold for any $\xi\in T_{p}M$
(and any $p\in M$)
\[
x^{a}\left(  \xi\right)  =x^{a}\left(  p\right)  \text{, }\xi=\sum
u^{a}\left(  \xi\right)  E_{a}\left(  p\right)
\]

Also we have the induced $\left(  x^{a},\overset{.}{x}^{a}\right)  $ pure
coordinates on $TM$ with
\[
\xi=\sum\overset{.}{x}^{a}\left(  \xi\right)  \partial_{x^{a}}%
\]

Let $\pi:TM\rightarrow M$ be the canonical projection given locally by
$\left(  x^{a},u^{a}\right)  \rightarrow\left(  x^{a}\right)  $.

The geodesic spray is the vectorfield $\Gamma$ on $TM$ whose integral curves
project down to the geodesics of $M$. Using mixed coordinates we may write:%
\[
\Gamma=\sum\overset{.}{x}^{a}\partial_{x^{a}}-\sum\Gamma_{ab}^{c}u^{a}%
u^{b}\partial_{u^{c}}%
\]
The projection on $M$ of the integral curves of $S=\tau\Gamma$ are
pregeodesics and we get:%
\begin{align*}
S  &  =\left(  \tau\overset{.}{x}^{a}\right)  \partial_{x^{a}}-\left(
\tau\Gamma_{ij}^{c}\right)  u^{i}u^{j}\partial_{u^{c}}-\left(  \tau\Gamma
_{mj}^{a}+\tau\Gamma_{jm}^{a}\right)  u^{m}u^{j}\partial_{u^{a}}\\
&  -\tau\Gamma_{mm}^{k}\left(  u^{m}\right)  ^{2}\partial_{u^{k}}-\tau
\Gamma_{mm}^{m}\left(  u^{m}\right)  ^{2}\partial_{u^{m}}%
\end{align*}
and using (\ref{Gammas}) we obtain%
\[
\left.  S\left(  E_{m}\right)  \right\vert _{D^{\infty}}=\frac{1}{2}\left.
E_{m}\left(  \tau\right)  \right\vert _{D^{\infty}}%
\]

Let $h$ be:%
\begin{equation}
h=E_{m}\left(  \tau\right)  /2\text{ (suppose for example }h<0\text{)}
\label{h}%
\end{equation}
we consider now the vectorfields%
\[%
\begin{tabular}
[c]{l}%
$H=hu^{a}\partial_{u^{a}}$\\
$A=\overset{.}{x}^{a}\partial_{x^{a}}-\sum_{\left(  a,b,c\right)  \neq\left(
m,m,m\right)  }\Gamma_{ab}^{c}u^{a}u^{b}\partial_{u^{c}}$\\
$B=h\left(  u^{m}\right)  ^{2}\partial_{u^{m}}$%
\end{tabular}
\ \
\]
and construct
\[
\widetilde{S}=S-H=\tau A+B-H
\]
and still the integral curves of $\widetilde{S}$ project down on pregeodesics
in $M$.

Now fix $p\in D^{\infty}$ Then $\xi=E_{m}\left(  p\right)  $ is stationary
point of $\widetilde{S}$ since
\begin{equation}
\tau\left(  p\right)  =0,u^{i}\left(  \xi\right)  =0,u^{m}\left(  \xi\right)
=1 \label{Psi}%
\end{equation}
and
\[
\widetilde{S}\left(  \xi\right)  =\tau\left(  p\right)  A\left(  \xi\right)
+B\left(  \xi\right)  -H\left(  \xi\right)  =0+h\left(  p\right)
\partial_{u^{m}}-h\left(  p\right)  \partial_{u^{m}}=0
\]

We linearize $\widetilde{S}$ at $\xi$ to obtain $\left.  D\widetilde
{S}\right\vert _{\xi}:T_{\xi}TM\rightarrow T_{\xi}TM$. Since $D\widetilde
{S}=d\tau\otimes A+\tau DA+DB-DH$ and taking account (\ref{Psi}) and
\[%
\begin{tabular}
[c]{c}%
$DB=dh\otimes\left(  u^{m}\right)  ^{2}\partial_{u^{m}}+2hu^{m}du^{m}%
\otimes\partial_{u^{m}}$\\
$DH=dh\otimes u^{a}\partial_{u^{a}}+h\otimes du^{a}\otimes\partial_{u^{a}}$%
\end{tabular}
\]
we have
\begin{align*}
\left.  D\widetilde{S}\right\vert _{\xi}  &  =\left(  \left.  d\tau\right\vert
_{p}\circ\pi_{\ast}\right)  \otimes\left(  \xi-\sum\Gamma_{mm}^{k}%
\partial_{u^{k}}\right) \\
&  +h\left(  p\right)  du^{m}\otimes\partial_{u^{m}}-h\left(  p\right)
du^{i}\otimes\partial_{u^{i}}%
\end{align*}

We remark that we have identify $\xi\in TM$ with $\xi\in T_{\xi}TM$ through
the $\left(  x^{a}\right)  $-canonical immersion $T_{p}M\hookrightarrow
T_{\xi}TM$ such that $\left.  \partial_{x^{a}}\right\vert _{p}\rightarrow
\left.  \partial_{x^{a}}\right\vert _{\xi}$. We compute now the eigenspaces to
apply later an stable manifold theorem.

\begin{itemize}
\item If $\eta\in T_{p}D^{\infty}\hookrightarrow T_{\xi}TM$ then
$0=d\tau\left(  \eta\right)  =u^{a}\left(  \eta\right)  $ and $\left.
D\widetilde{S}\right\vert _{\xi}\left(  \eta\right)  =0.$

\item If $\eta=\left.  \partial_{u^{i}}\right\vert _{\xi}$ then $0=dx^{i}%
\left(  \eta\right)  =d\tau\left(  \eta\right)  =du^{a}\left(  \eta\right)  $
($a\neq i)$ and $du^{i}\left(  \eta\right)  =1$ thus $\left.  D\widetilde
{S}\right\vert _{\xi}\left(  \eta\right)  =-h\left(  p\right)  \eta$

\item If $\eta=\left.  \partial_{u^{m}}\right\vert _{\xi}$ then $\left.
D\widetilde{S}\right\vert _{\xi}\left(  \eta\right)  =h\left(  p\right)  \eta$

\item $\left.  D\widetilde{S}\right\vert _{\xi}\left(  \xi\right)  =2h\left(
p\right)  \left(  \xi-\sum\Gamma_{mm}^{k}\left(  p\right)  \partial_{u^{k}%
}\right)  $. It is easy to prove that there exist an eigenvector associated to
eigenvalue $2h\left(  p\right)  $ with the look $\eta=\xi-\sum c_{i}%
\partial_{u^{i}}$ for some constant $c_{i}$.
\end{itemize}

Collected the information we have eigenvalues $0,-h\left(  p\right)  ,$
$h\left(  p\right)  $, $2h\left(  p\right)  $ with multiplicity $m-1,$ $m-1$,
$1$ $,1$respectively. (their are $2m$ counting their multiplicities). Because
the eigenvalue $2h\left(  p\right)  $ is smaller than (and not equal to) any
other negative eigenvalue we conclude (by certain refinement of the stable
manifold theorem) that there exist a $\widetilde{S}$-stable line
$\widetilde{L}\subset TM$ with $\xi\in\widetilde{L}$ and $T_{\xi}\widetilde
{L}=Span\left(  \xi-\sum c_{i}\partial_{u^{i}}\right)  $. The projection
$L=\pi\left(  \widetilde{L}\right)  $ of this stable line sweeps out the
smoothly immersed pregeodesic.

\subsection{The natural equation for $D^{\infty}$.\label{Sec32}}

We find a \ (locally) coordinate system $\left(  z^{a}\right)  $ such that
$\left(  g_{ab}\right)  =\left(  g\left(  \partial_{z^{a}},\partial_{z^{b}%
}\right)  \right)  $ is as%
\begin{equation}
\left(  g_{ab}\right)  =\left(
\begin{array}
[c]{cc}%
\left(  g_{ij}\right)  & 0\\
0 & 1/z^{m}%
\end{array}
\right)  \label{gij}%
\end{equation}

Note that in these coordinates $\partial_{z^{m}}$ be a polar normal
vectorfield and their integral curves $\overline{\gamma}=\overline{\gamma
}\left(  s\right)  :\left\{  z^{i}=cte\text{, }z^{m}=s\right\}  $ are the
pregeodesics of the previous section. Thus the $z^{m}\left(  =s\right)
$-coordinate parametrize the pregeodesic $\gamma$ in such way that%
\begin{equation}
\left.  g\left(  \partial_{z^{m}},\partial_{z^{m}}\right)  \right\vert
_{\gamma}=g\left(  \frac{d\overline{\gamma}}{ds},\frac{d\overline{\gamma}}%
{ds}\right)  =\frac{1}{s} \label{dgamma}%
\end{equation}

Therefore as previous question, we analyze the existence of such (canonical)
parametrization. Next we will construct from the parameter $s$ (of obvious
way) a natural equation $(z^{m}=0)$ for $D^{\infty}$.

We start with a fixed parametrization on the polar $D^{\infty}-$transversal
pregeodesic $\gamma=\gamma\left(  t\right)  $ defined for $\left\vert
t\right\vert <\varepsilon$, with $\gamma\left(  0\right)  =p\in D^{\infty}$.
Using by example a polar normal vectorfield which has $\gamma$ as integral
curve, is easy to see that the function
\[
\Phi\left(  t\right)  =\frac{1}{g\left(  \gamma^{\prime}\left(  t\right)
,\gamma^{\prime}\left(  t\right)  \right)  }\text{, for }0<\left\vert
t\right\vert <\varepsilon\text{, }\Phi\left(  0\right)  =0
\]
is a $C^{\infty}$-function on the whole interval $\left(  -\varepsilon
,\varepsilon\right)  $ .with $\Phi^{\prime}\left(  0\right)  \neq0$ and we may
write
\[
\Phi\left(  t\right)  =t\Psi\left(  t\right)
\]
for some $C^{\infty}$-function $\Psi$ such that $\Psi\left(  t\right)  \neq0$
for $\left\vert t\right\vert <\varepsilon$ (suppose for example $\Psi\left(
t\right)  >0$) . Let
\begin{equation}
\psi\left(  t\right)  =\frac{1}{2\sqrt{\Psi\left(  t\right)  }} \label{psi}%
\end{equation}

We find a new parameter $s=\mathbf{s}\left(  t\right)  $, with inverse
$t=\mathbf{t}\left(  s\right)  $, such that the reparametrized curve
$\overline{\gamma}(s)=\gamma\left(  \mathbf{t}\left(  s\right)  \right)  $
verify (\ref{dgamma}). Thus%
\[
\frac{1}{s}=g\left(  \frac{d\overline{\gamma}}{ds},\frac{d\overline{\gamma}%
}{ds}\right)  =\left(  \frac{d\mathbf{t}}{ds}\right)  ^{2}g\left(
\frac{d\gamma}{dt},\frac{d\gamma}{dt}\right)
\]
and therefore the function $s=\mathbf{s}\left(  t\right)  $ satisfy the
differential equation of separate variables
\[
\frac{ds}{\sqrt{s}}=\frac{\psi dt}{2\sqrt{t}}%
\]
integrating both members we have an explicit solution:%
\begin{equation}
\mathbf{s}\left(  t\right)  =sgn\left(  t\right)  \left(  \int_{0}^{t}%
\frac{\psi dx}{\sqrt{x}}\right)  ^{2} \label{s}%
\end{equation}

\begin{lemma}
\label{Lemma}Let $\psi:\left(  -\varepsilon,\varepsilon\right)  \rightarrow
\mathbb{R}$ be a $C^{\infty}$-function. Then $\mathbf{s}\left(  t\right)  $
defined in $\left(  -\varepsilon,\varepsilon\right)  $ as (\ref{s}) as well is
a $C^{\infty}$-function
\end{lemma}

\begin{proof}
See appendix
\end{proof}

Now there exist an unique smooth function $z^{m}$ defined over a neighborhood
of $D^{\infty}$ such that $z^{m}\left(  \overline{\gamma}\left(  s\right)
\right)  =s$, where $\overline{\gamma}=\overline{\gamma}\left(  s\right)  $ is
any normal pregeodesic and $s$ is their natural parameter.

\begin{remark}
\label{Rem2b}Using the definition of $\psi$ in (\ref{psi}) we have
\[
\frac{\psi\left(  t\right)  }{\sqrt{t}}=\pm\frac{1}{2}\sqrt{\left\vert
g\left(  \gamma^{\prime}\left(  t\right)  ,\gamma^{\prime}\left(  t\right)
\right)  \right\vert }%
\]
and for $x=\gamma\left(  t_{0}\right)  $ belonging to such neighborhood
$4z^{m}\left(  x\right)  =\mathbf{s}\left(  t_{0}\right)  $ is the signed
square of the the arc length (if $t_{0}>0)$, or of the proper time (if
$t_{0}\leq0$) to the normal pregeodesic segment between $x$ and $D^{\infty}$
\end{remark}

\subsection{The natural polar coordinates.\label{Sec33}}

Rewriting with more formalism the end of the previous section we may say that:

For some $\varepsilon>0$ there exist a smooth function $\zeta:D^{\infty}%
\times\left(  -\varepsilon,\varepsilon\right)  \rightarrow M$, such that
$s\rightarrow\zeta\left(  x,s\right)  $ ($\left\vert s\right\vert
<\varepsilon$) is the normal pregeodesic at $x=\zeta\left(  x,0\right)  \in
D^{\infty}$ with the natural parametrization. Since $\zeta$ is non singular at
the points $\left(  x,0\right)  $, we may suppose (replacing $M$ by some
neighborhood of $D^{\infty})$ that $\zeta$ is diffeomorphism. We write the
inverse as $\zeta^{-1}=\left(  \sigma,z^{m}\right)  :M\rightarrow D^{\infty
}\times\left(  -\varepsilon,\varepsilon\right)  $, and the vectorfield
$\partial_{z^{m}}$ defined by%
\[
\left.  \partial_{z^{m}}\right\vert _{p}=\left.  \frac{\partial\zeta}{\partial
s}\right\vert _{\left(  \sigma\left(  p\right)  ,z^{m}\left(  p\right)
\right)  }%
\]
is called polar-normal pregeodesic vectorfield. Note that using an auxiliary
polar normal frame $\left(  E_{i},\partial_{z^{m}}\right)  $ and the last
equation of (\ref{Gammas}) (here is $\tau=z^{m}$) we see that
\begin{equation}
\nabla_{\partial_{z^{m}}}\partial_{z^{m}}=-\frac{1}{2z^{m}}\partial_{z^{m}}
\label{zmzm}%
\end{equation}
Moreover by (\ref{dgamma}) we have%
\begin{equation}
g\left(  \partial_{z^{m}},\partial_{z^{m}}\right)  =\frac{1}{z^{m}}
\label{gmm}%
\end{equation}

We start now with a coordinate system $\left(  x^{i}\right)  $ on $D^{\infty}$
We will prove that the coordinates $\left(  z^{i},z^{m}\right)  $ on $M$ where
$z^{i}=x^{i}\circ\sigma$, are natural coordinates, that is $\left(
g_{ab}\right)  $ is as (\ref{gij}). Taking account (\ref{gmm}) we have%
\[
0=\frac{\partial}{\partial z^{i}}\left(  \frac{1}{z^{m}}\right)
=\frac{\partial}{\partial z^{i}}g\left(  \partial_{z^{m}},\partial_{z^{m}%
}\right)  =2\square_{\partial_{z^{i}}}\partial_{z^{m}}\left(  \partial_{z^{m}%
}\right)
\]
and therefore, using also (\ref{zmzm})
\begin{align*}
\frac{\partial g_{im}}{\partial z^{m}}  &  =\frac{\partial}{\partial z^{m}%
}g\left(  \partial_{z^{i}},\partial_{z^{m}}\right)  =g\left(  \nabla
_{\partial_{z^{m}}}\partial_{z^{m}},\partial_{z^{i}}\right) \\
&  =-\frac{1}{2z^{m}}g_{im}%
\end{align*}
By the Corollary \ref{Cor1}, $g_{im}$ is a differentiable function on $M$. In
order to prove that $g_{im}=0$, we fix the variables $z^{i}$. The smooth
function $\phi(t)=g_{im}\left(  z^{i},t\right)  $ defined on some open
interval around zero satisfy.
\begin{equation}
-2t\frac{d\phi}{dt}=\phi\text{ for all }t \label{dfidt}%
\end{equation}
the following Lemma proves that but always that $\phi$ is identically null

\begin{lemma}
Let $\phi:I\rightarrow\mathbb{R}$ a smooth function defined over some open
interval $I$ around zero which satisfy (\ref{dfidt}). Then $\phi(t)=0$
$\forall t$.
\end{lemma}

\begin{proof}
First note that by (\ref{dfidt}) is $\phi(0)=0$. Suppose that there exist
$t_{0}\in I$ with $\phi(t_{0})\neq0$ (for example $\phi(t_{0})>0$.$)$ Let
$J=\left(  a,b\right)  \subset I\cap\mathbb{R}^{+}$ an open maximal interval
containing $t_{0}$ such that $\phi(t)>0$for all $t\in J$. Then using the
differential equation (\ref{dfidt}) it is easy to see that there exist a
constant $C>0$ such that $\phi(t)=C/\sqrt{t}$. Then $a>0$ (if not
$\lim_{t\rightarrow0^{+}}\phi(t)=+\infty$ $\neq0=$ $\phi(0)$), but then
$\phi(a)=C/\sqrt{a}>0$ and thus $J$ is not maximal.
\end{proof}

Therefore we have established \ the existence of polar-normal coordinates as
explain in the following

\begin{theorem}
Around each point of $D^{\infty}$ there exist a coordinate system $\left(
z^{i},z^{m}\right)  $ such that $\left(  g_{ij}\right)  =\left(  \left(
\partial_{z^{i}},\partial_{z^{j}}\right)  \right)  $ is as (\ref{gij})
\end{theorem}

\section{The Curvature near to $D^{\infty}$.\label{Sec4}}

In order to analyze the limiting behavior of the curvatures on $M-D^{\infty}$
as we approach the polar hypersurface $D^{\infty}$, we fix a polar-normal
coordinate system $\left(  z^{i},z^{m}=\tau\right)  $ whose domain is the
whole space $M$ (if not we restrict to the domain). Let $\left(
g_{ab}\right)  $ be as (\ref{gij}) with $z^{m}=\tau$. The inverse is
\[
\left(  g_{ab}\right)  ^{-1}=\left(  g^{ab}\right)  =\left(
\begin{array}
[c]{cc}%
\left(  g^{ij}\right)  & 0\\
0 & \tau
\end{array}
\right)
\]

Recall that here the Christopher symbols are
\[%
\begin{tabular}
[c]{l}%
$\Gamma_{cab}=\left(  \square_{\partial_{z^{a}}}\partial_{z^{b}}\right)
\left(  \partial_{z^{c}}\right)  =\dfrac{1}{2}\left\{  \dfrac{\partial g_{ac}%
}{\partial z^{b}}+\dfrac{\partial g_{bc}}{\partial z^{a}}-\dfrac{\partial
g_{ab}}{\partial z^{c}}\right\}  $\\
$\Gamma_{ab}^{c}=\sum\Gamma_{dab}g^{dc}$%
\end{tabular}
\
\]

The contravariant curvature is
\[
R\left(  A,B\right)  C=\nabla_{A}\nabla_{B}C-\nabla_{B}\nabla_{A}%
C-\nabla_{\left[  A,B\right]  }C
\]
their components $R_{cab}^{d}$ are defined by $R\left(  \partial_{z^{a}%
},\partial_{z^{b}}\right)  \partial_{z^{c}}=\sum R_{cab}^{d}\partial_{z^{d}}$
and are given by%
\[
R_{cab}^{d}=\frac{\partial\Gamma_{bc}^{d}}{\partial z^{a}}-\frac
{\partial\Gamma_{ac}^{d}}{\partial z^{a}}+\sum\Gamma_{bc}^{e}\Gamma_{ae}%
^{d}-\sum\Gamma_{ac}^{e}\Gamma_{be}^{d}%
\]

The covariant Ricci tensor is $Ric\left(  A,B\right)  =tr\left\{  V\rightarrow
R\left(  A,V\right)  CB\right\}  $ which has components
\[
R_{ca}=\sum R_{cad}^{d}\text{, }R_{c}^{d}=\sum R_{ca}g^{ad}%
\]
where $R_{c}^{d}$ are the components of the contravariant Ricci tensor defined
by the identity $Ric\left(  A,B\right)  =g\left(  RIC\left(  A\right)
,B\right)  $.

Finally the components of the Weil curvature tensor $W$ are:%
\begin{align*}
W_{cab}^{d}  &  =R_{cab}^{d}+\frac{1}{m-2}\left(  \delta_{a}^{d}R_{cb}%
-\delta_{b}^{d}bR_{ca}+g_{cb}R_{a}^{d}-g_{ca}R_{b}^{d}\right) \\
&  +\frac{S}{\left(  m-1\right)  \left(  m-2\right)  }\left(  \delta_{b}%
^{d}g_{ca}-\delta_{a}^{d}g_{cb}\right)
\end{align*}
\bigskip

where $S$ is the scalar curvature.
\[
S=tr\left(  RIC\right)  =\sum R_{c}^{c}%
\]
By inspection of previous formulas we have the following:

\begin{lemma}
In the local natural coordinates we have:

\begin{enumerate}
\item $\Gamma_{cab}\cong0$ except $\Gamma_{mmm}=-1/\tau^{2}$

\item $\Gamma_{ab}^{c}\cong0$ except $\Gamma_{mm}^{m}=-1/\tau$

\item $R_{abc}^{d}\cong0$, except perhaps $-R_{mmj}^{i}=R_{mjm}^{i}%
\cong-\Gamma_{jm}^{i}/\tau$

\item $R_{ab}\cong0,$except perhaps $R_{mm}=\sum R_{mmj}^{j}$, but $\tau
R_{mm}\cong0$

\item $R_{a}^{b}\cong0$, $S\cong0$

\item $W_{abc}^{d}\cong0$ except perhaps $W_{mjm}^{d}$ but $\tau W_{mjm}%
^{d}\cong0$
\end{enumerate}
\end{lemma}

\bigskip

As consequence we have the following

\begin{theorem}
\label{Th2}Let $R$, $Ric$, $RIC$, $W$ be the previous curvature tensors of
$\left(  M-D^{\infty},g\right)  $, and $X,Y,Z\in\mathfrak{X}\left(  M\right)
$, and let $\tau$ be the natural equation of $D^{\infty}$. Denote by
$R^{\infty}$,...etc. the corresponding tensors to the Riemannian space
$D^{\infty}$

\begin{enumerate}
\item $\tau R$, $\tau Ric$, $RIC$, $\tau W$ are $C^{\infty}-$tensorfields
defined around $D^{\infty}$

\item $R\left(  X,Y\right)  Z\cong0$ if $Z\in\mathfrak{X}_{M}\left(
D^{\infty}\right)  $, or $X,Y\in\mathfrak{X}_{M}\left(  D^{\infty}\right)  $.
Moreover if $X,Y,Z\in\mathfrak{X}_{D^{\infty}}\left(  M\right)  $ then
$\left.  R\left(  X,Y\right)  Z\right\vert _{D^{\infty}}=R^{\infty}\left(
X,Y\right)  Z$

\item $Ric\left(  X,Y\right)  \cong0$ if $X$ or $Y$ are tangent to $D^{\infty
}$.

\item $W\left(  X,Y\right)  Z\cong0$ if $Z\in\mathfrak{X}_{M}\left(
D^{\infty}\right)  $, or $X,Y\in\mathfrak{X}_{D^{\infty}}\left(  M\right)  $
\end{enumerate}
\end{theorem}

\begin{proof}
It is straightforward. For example we will prove 2:

Note that $R\left(  \partial_{z^{a}},\partial_{z^{b}}\right)  \partial_{z^{c}%
}=\sum R_{cab}^{d}\partial_{z^{d}}\cong0$ if $c\neq m$ or $c=m$ , $a\neq
m$,and $n\neq m$. Writing $X=\sum X^{a}\partial_{z^{a}}$ ...etc. then%
\begin{align*}
R\left(  X,Y\right)  Z  &  =X^{a}Y^{b}Z^{i}R\left(  \partial_{z^{a}}%
,\partial_{z^{b}}\right)  \partial_{z^{i}}+X^{i}Y^{j}Z^{m}R\left(
\partial_{z^{i}},\partial_{z^{j}}\right)  \partial_{z^{m}}\\
&  \cong X^{m}Y^{j}Z^{m}R\left(  \partial_{z^{m}},\partial_{z^{j}}\right)
\partial_{z^{m}}+X^{i}Y^{m}Z^{m}R\left(  \partial_{z^{i}},\partial_{z^{m}%
}\right)  \partial_{z^{m}}%
\end{align*}

But if $Z$ is tangent to $D^{\infty}$ is $Z^{m}=\tau C^{m}$ for some smooth
$C^{m}$...etc. The proof finish taking account that $\tau R\cong0$
\end{proof}

\begin{remark}
It is easy to see that on the natural coordinates $\left(  z^{i},\tau\right)
$, the condition to $R\cong0$ is equivalent to
\[
\left.  \frac{\partial g_{ij}}{\partial z^{m}}\right\vert _{z^{m}=0}=0
\]
and this condition may be surely expressed without coordinates.\newpage
\end{remark}

\section{Conformal Geometry.}

We consider the conformal class $\left(  M,\mathcal{C}\right)  $ of a
Riemann-Lorentz space $\left(  M,g\right)  $ with polar end $D^{\infty}$. We
recall that
\[
\mathcal{C}=\left\{  e^{2\sigma}g:\sigma\in C^{\infty}\left(  M\right)
\right\}
\]
We remark that $\left(  M,\overline{g}\right)  $ is also a Riemann-Lorentz
space with polar end $D^{\infty}$ for any $\overline{g}\in\mathcal{C}$.

Of course the polar-normal pregeodesic introduced in the section \ref{Sec31}
are not determined by the conformal class $\mathcal{C}$. Also the polar-normal
direction on $D^{\infty}$ is not determined by $\mathcal{C}$. In fact, if
$\left(  z_{i},z_{m}\right)  $ are the canonical coordinates associated to
$g$, then $\overline{g}=fg$ with $f=f\left(  z_{i},z_{m}\right)  >0$, then
$\left(  \partial_{z_{i}},\partial_{z_{m}}\right)  $ is still a polar (not
orthonormal) frame and (see subsection \ref{Sec22})%
\[
\overline{\Gamma}_{kmm}=\frac{1}{2}\frac{\partial\left(  f/z_{m}\right)
}{\partial z_{k}}=\frac{1}{2z_{m}}\frac{\partial f}{\partial z_{k}}%
\]
therefore $\left.  z_{m}\overline{\Gamma}_{kmm}\right\vert _{z_{m}=0}$ may to
take any arbitrary value moving $f$. Of course the polar normal direction
remains invariant if $\left.  \dfrac{\partial f}{\partial z_{k}}\right\vert
_{z_{m}=0}=0$.

However we will prove that the family of the polar-normal pregeodesics are the
same for all the metrics of the conformal subclass
\begin{equation}
\mathcal{C}_{g}^{\prime}=\left\{  e^{2\sigma}g:\sigma=f\circ\tau_{g}\text{
with }f\in C^{\infty}\left(  \mathbb{R}\right)  \right\}  \label{Cprimag}%
\end{equation}
where $\tau_{g}=0$ is the canonical equation for $D^{\infty}$ (with respect to
$\left(  M,g\right)  $ $)$ established in section \ref{Sec32}$.$Moreover
$\mathcal{C}^{\prime}=\mathcal{C}_{g}^{\prime}$ depends only to the family of
the polar-normal pregeodesic and not of the initial metric $g$. This means
that: if $g,\overline{g}\in\mathcal{C}$ then $\overline{g}$ has the same polar
normal pregeodesics that $g$, if and only if $\overline{g}\in\mathcal{C}%
_{g}^{\prime}$.

On the other hand the family of hypersurfaces $D_{t}^{\infty}$ with equation
$\left(  \tau_{g}=t\right)  $ depends only to $\mathcal{C}_{g}^{\prime}$. This
is the family of the \emph{simultaneity hypersurfaces} which determines the
\emph{simultaneity distribution} $\mathcal{D}_{g}$. We define a (abstract)
simultaneity distribution as a completely integrable $\left(  m-1\right)
$-distribution $\mathcal{D}$ such that $D^{\infty}$ is an integral manifold. \

Let $N$ be an everywhere non isotropic vectorfield on $M$, transverse to
$D^{\infty}$. We define the distribution $N^{\bot}$ as $N^{\bot}\left(
p\right)  =N\left(  p\right)  ^{\bot}$ if $p\in M-D^{\infty}$, and $N^{\bot
}\left(  p\right)  =T_{p}D^{\infty}$ if $p\in D^{\infty}$.

\begin{lemma}
\label{Lem2}If $N$ is an everywhere non isotropic vectorfield on $M$
transverse to $D^{\infty}$, then $N^{\bot}$ is a smooth distribution, and for
any $g\in\mathcal{C}$, the function $g\left(  N,N\right)  ^{-1}$ extend to
$D^{\infty}$, and $g\left(  N,N\right)  ^{-1}=0$ is an equation for
$D^{\infty}$. Reciprocally if $\mathcal{D}$ is a simultaneity distribution on
the conformal space $\left(  M,\mathcal{C}\right)  $, then there exist $N$ non
isotropic vectorfield on $M$ transverse to $D^{\infty}$such that
$\mathcal{D}=N^{\bot}$.
\end{lemma}

\begin{proof}
Fix any $g\in\mathcal{C}$ and let $\left(  z_{i},z_{m}\right)  $ be polar
normal coordinates around $D^{\infty}$ this means that%
\[
\left(  g_{ab}\right)  =\left(
\begin{array}
[c]{cc}%
\left(  g_{ij}\right)  & 0\\
0 & 1/z_{m}%
\end{array}
\right)
\]
suppose
\begin{equation}
N=\sum\lambda_{i}\partial_{z_{i}}+\lambda_{m}\partial_{z_{m}} \label{N}%
\end{equation}
The non isotropic condition for $N$ assure that $N^{\bot}$ is a $\left(
m-1\right)  -$distribution away $D^{\infty}$. We may describe $N^{\bot}$ on
$M-D^{\infty}$ as the family of vectorfields $X=\sum X_{i}\partial_{z_{i}%
}+X_{m}\partial_{z_{m}}$, such that
\begin{equation}
\sum_{j}\Lambda_{j}X_{j}+\lambda_{m}X_{m}=0\text{ with }\Lambda_{j}=\sum
_{i}z_{m}g_{ij}\lambda_{i} \label{Landas}%
\end{equation}
Note that the coefficients $\Lambda_{j}$ and $\lambda_{m}$ are smooth on $M$.
But for $z_{m}=0$ are $\Lambda_{j}\left(  z_{i},0\right)  =0$, and
$\lambda_{m}\left(  z_{i},0\right)  \neq0$ (by transversality). Previous
equation gives $X_{m}=0.$ This means that the vectorfield $X$ is tangent to
$D^{\infty}$ and that the condition (\ref{Landas}) define the whole
distribution $N^{\bot}$. This proves that $N^{\bot}$ is smooth. Also
\[
\frac{1}{g\left(  N,N\right)  }=\frac{z_{m}}{z_{m}\Lambda+\lambda_{m}^{2}%
}\text{ with }\Lambda=\sum_{i}g_{ij}\lambda_{i}\lambda_{j}%
\]
since $\lambda_{m}\left(  z_{i},0\right)  \neq0$, this proves that $g\left(
N,N\right)  ^{-1}$ extend to $D^{\infty}$, and $g\left(  N,N\right)  ^{-1}=0$
is an equation for $D^{\infty}.$

We may describe a simultaneity distribution $\mathcal{D}$ in the polar normal
coordinates $\left(  z_{i},z_{m}\right)  $ as the vector fields $X=\sum
X_{i}\partial_{z_{i}}+X_{m}\partial_{z_{m}}$ which satisfy a condition as%
\begin{equation}
\sum_{j}\Lambda_{j}X_{j}+\lambda_{m}X_{m}=0 \label{Landas1}%
\end{equation}

Since $D^{\infty}$ is integral manifold we conclude that $\Lambda_{j}\left(
z_{i},0\right)  =0$ and there exist smooth functions $\mu_{j}$ such that
\[
\Lambda_{j}=z_{m}\mu_{j}%
\]
since $\left(  g_{ij}\right)  $ is nonsingular there exist $\lambda_{i}$ such
that $\mu_{j}=\sum g_{ij}\lambda_{i}$, and (\ref{Landas1}) becomes in
(\ref{Landas}), select $N$ as in (\ref{N}) we conclude $\mathcal{D}=N^{\bot}$.
\end{proof}

\begin{theorem}
\label{Th3}Given an abstract simultaneity distribution $\mathcal{D}$\ then
there exist $\overline{g}\in$ $\mathcal{C}$ such that $\mathcal{D=D}%
_{\overline{g}}$.
\end{theorem}

\begin{proof}
We take an auxiliary metric $g\in\mathcal{C}$. By previous Lemma we may select
an everywhere non null vectorfield $N$ such that $g\left(  N,N\right)
^{-1}=0$ is an equation for $D^{\infty}$, and $N^{\bot}=\mathcal{D}$. We fix a
point $x_{0}\in D^{\infty}$, and let $\gamma_{x_{0}}:\left(  -c,c\right)
\rightarrow M$ \emph{any} regular parametrization of the integral curve
$\alpha_{x_{0}}$ of $N$ which $\alpha_{x_{0}}\left(  0\right)  =x_{0}%
=\gamma_{x_{0}}\left(  0\right)  $ (we may take for example $\gamma_{x_{0}%
}=\alpha_{x_{0}}$). Let $\Sigma_{t}$ the integral manifold of $\mathcal{D}$ by
$\gamma_{x_{0}}\left(  t\right)  $. For any $x\in\Sigma$ we parametrize the
integral curve $\alpha_{x}$ of $N$ by $\gamma_{x}:\left(  -c,c\right)
\rightarrow M$ such that $\gamma_{x_{0}}\left(  t\right)  \in\Sigma_{t}$.
\ Let $\Phi:\Sigma\times\left(  -c,c\right)  \rightarrow M$ be such that
$\Phi\left(  x,t\right)  =\gamma_{x}\left(  t\right)  $. It is straightforward
to see that $\Phi$ is not singular over points of $\Sigma\times\left\{
0\right\}  $, and we may suppose without lost generality that $\Phi$ is
diffeomorphism. Using the inverse $\Phi^{-1}:M\rightarrow D^{\infty}%
\times\left(  -c,c\right)  $ we may make the unique coordinate system $\left(
x_{i},x_{m}\right)  $ on $M$, such that%
\[
\Phi:\left\{
\begin{array}
[c]{c}%
x_{i}=u_{i}\\
x_{m}=t
\end{array}
\right.
\]
are the equations of $\Phi$ ( here $\left(  u_{i}\right)  $ are a fixed
coordinate system on $D^{\infty}$ and $t$ is the coordinate on $\left(
-c,c\right)  $). Since the curves $\left\{  x_{i}=cte\text{, }x_{m}=t\right\}
$ are preintegral curves of $N$ we conclude that $\partial_{x_{m}}%
=e^{-\varphi}N$ for some smooth $\varphi$. Also since (for a fixed $t$)
$x_{m}=t$ is the equation of $\Sigma_{t}$ then $\partial_{x_{i}}%
\in\mathcal{D=}\partial_{x_{m}}^{\bot}$ and the metric $g$ in coordinates
$\left(  x_{i},x_{m}\right)  $ has the matrix
\[
\left(  g_{ab}\right)  =\left(
\begin{array}
[c]{cc}%
\left(  g_{ij}\right)  & 0\\
0 & e^{-2\varphi}g\left(  N,N\right)
\end{array}
\right)
\]

Since $g\left(  N,N\right)  ^{-1}=0$ is an equation for $D^{\infty}$, as
$x_{m}=0$ we conclude that $g\left(  N,N\right)  =h/x_{m}$ where $h>0$
everywhere. then the matrix of $\overline{g}=e^{2\varphi}h^{-1}g$ with respect
to such coordinates are
\[
\left(  \overline{g}_{ab}\right)  =\left(
\begin{array}
[c]{cc}%
\left(  e^{2\varphi}g_{ij}\right)  & 0\\
0 & 1/x_{m}%
\end{array}
\right)
\]
and we conclude that $\left(  x_{i},x_{m}\right)  $ are polar normal
coordinates for $\overline{g}$, and $\mathcal{D}=\mathcal{D}_{\overline{g}}$.
\end{proof}

\begin{remark}
\label{Rem5} With the hypothesis of previous theorem, the same argument proves
that the class $\mathcal{C}_{\mathcal{D}}$ of all $\overline{g}\in\mathcal{C}$
such that $\mathcal{D=D}_{\overline{g}}$ it is equal to $\mathcal{C}%
_{g}^{\prime}.$ The key is that we are free to parametrize $\gamma_{x_{0}}$.
This means that the previous coordinate $x_{m}$ (and therefore $\overline{g}$)
is determined up composition by arbitrary diffeomorphism. $f$ $\in C^{\infty
}\left(  \mathbb{R}\right)  $
\end{remark}

\subsection{Cosmological remarks.}

We consider the conformal class $\left(  M,\mathcal{C}\right)  $ of a
Riemann-Lorentz space $\left(  M,g\right)  $ with polar end $D^{\infty}$. This
is the support to a causality structure of the Lorentz component $D^{-}$. The
aim of this section, is to know if it is possible to find a
big-bang.\emph{cosmologically} privileged metric (around $D^{\infty}$). Of
Course we must to impose to $\left(  M,\mathcal{C}\right)  $ some initial
restriction \ as for example that should be $D^{\infty}$ conformal flat.

We recall that a Robertson-Walker space is a warped product $I\times
_{f}S=\left(  I\times S,g_{RW}\right)  $ where $I=\left(  0,t^{\ast}\right)  $
is an open interval and $\ f:I\rightarrow\mathbb{R}$ is a smooth function and
$\left(  S,g_{S}\right)  $ is a Riemannian manifold with constant sectional
curvature.$C_{0}$ Finally
\[
g_{RW}=-dt^{2}+f\left(  t\right)  ^{2}g_{S}=f\left(  t\right)  ^{2}g
\]

where $g=-f\left(  t\right)  ^{-2}dt^{2}+g_{S}$. Recall that $\left\{
t\right\}  \times S$ are simultaneity hypersurfaces of constant curvature
$C\left(  t\right)  =C_{0}f\left(  t\right)  ^{-2}$.

We remark that the flow $\zeta_{t}:S\rightarrow\left\{  t\right\}  \times S$ ,
$x\rightarrow\left(  t,x\right)  $ are homoteties of ratio $f\left(  t\right)
^{2}$

This suggest that in $\left(  M,g\right)  $, near to $D^{\infty}$ the metric
$g_{c}=-\left(  \tau_{g}\right)  g$ is the cosmologically relevant one. In
fact we have:

\begin{proposition}
\label{p4}Let $\left(  M,g\right)  $ be a four dimensional Riemann-Lorentz
space with polar end $D^{\infty}$ the flow $\zeta^{g}:D^{\infty}\times\left(
-\varepsilon,\varepsilon\right)  \rightarrow M$ in section \ref{Sec33} moves
$D^{\infty}$ by%
\[
\zeta_{t}^{g}:D^{\infty}\rightarrow D_{t}^{\infty}=\zeta^{g}\left(  D^{\infty
}\times\left\{  t\right\}  \right)
\]
Suppose that $D^{\infty}$ has constant (Riemannian) curvature and $\zeta
_{t}^{g}:D^{\infty}\rightarrow D_{t}^{\infty}$ are homoteties. Then the
simultaneity distribution $\mathcal{D}_{g}$ has integral manifolds which are
of constant curvature\footnote{Note that such property depends only of the
conformal subclass $\mathcal{C}_{g}^{\prime}$.} , and $g_{c}=-\left(  \tau
_{g}\right)  g$ becomes (locally) $D^{-}$ into a Robertson Walker space.
\end{proposition}

Note that $g_{c}$ induces the same causality structure as $\left[  g\right]  $
on $D^{-}$. However $g_{c}\notin\mathcal{C}_{g}$ since $\tau_{g}$ is null over
$D^{\infty}$.

Of course the physical relevant metric $g_{c}=-\left(  \tau_{g}\right)  g$ is
not determined by the conformal structure $\mathcal{C}=\mathcal{C}_{g}$. but
neither by the restricted conformal class $\mathcal{C}_{g}^{\prime}.$ However
the physical relevant metric are determined up the selection of an
\emph{universal time }$\tau$ in $\left(  M,\mathcal{C}\right)  $, as we was
proved in Remark \ref{Rem5}

\begin{remark}
\label{Rem6}To say that the integral hypersurfaces of a simultaneity
distribution $\mathcal{D}$ are of constant curvature has meaning into the
conformal Riemann Lorentz space space $\left(  M,\mathcal{C}\right)  $. This
means that their integral hypersurfaces has constant curvature with respect to
any (or some) $g$ belonging to the restricted conformal class $\mathcal{C}%
_{\mathcal{D}}^{\prime}$ induced by $\mathcal{D}$ according Theorem \ref{Th3}
\end{remark}

Finally we set out the following conjecture that explain the philosophical
motivation mentioned at the beginning of this subsection.

\begin{conjecture}
\label{p5}Let $\left(  M,\mathcal{C}\right)  $ be a Riemann-Lorentz conformal
space with polar end $D^{\infty}$. Suppose that $D^{\infty}$ is conformal-
flat. Then there exist (locally) a simultaneity distribution $\mathcal{D}$ of
constant curvature. Moreover $\mathcal{D}$ is univocally determined.
\end{conjecture}

Note that by Remark \ref{Rem5} we find (locally) a equation $\left(
\tau=0\right)  $ of $D^{\infty}$ whose level hypersurfaces $\left(
\tau=cte\right)  $ are of constant curvature (as we explain in the Remark
\ref{Rem6}). Then the conjecture says that the equation $\left(
\tau=0\right)  $ of $D^{\infty}$ is univocally determined by the constant
curvature condition up diffeomorphism. $\phi:\mathbb{R\rightarrow}\mathbb{R}$.

\section{Appendix}

This appendix is devoted to prove the following result:

\begin{lemma}
Let $\psi:I_{\varepsilon}=\left(  -\varepsilon,\varepsilon\right)
\rightarrow\mathbb{R}$ be a $C^{\infty}$-function. Then $F\left(  t\right)  $
defined in $I_{\varepsilon}$ as%
\[
F\left(  t\right)  =\epsilon\left(  t\right)  \left(  \int_{0}^{t}\frac{\psi
dx}{\sqrt{x}}\right)  ^{2}%
\]
is also a $C^{\infty}$-function,.where $\epsilon$ is the sign function
($\epsilon\left(  x\right)  =1$ if $x>0$, $\epsilon\left(  x\right)  =-1$ if
$x<0$, and $\epsilon\left(  0\right)  =0$)
\end{lemma}

\begin{proof}
To simplify we will denote the functions without reference to the variable,
$I_{\varepsilon}^{\ast}=I_{\varepsilon}-\left\{  0\right\}  $, $J$ denote the
absolute value function (that is $J\left(  x\right)  =\left\vert x\right\vert
=\epsilon\left(  x\right)  x$ and $\int g$ denote
\[
\left(  \int g\right)  \left(  t\right)  =\int_{0}^{t}g\text{ for }\left\vert
t\right\vert <\varepsilon\text{, and }g\text{ integrable in }\left(
-\varepsilon,\varepsilon\right)
\]

$g^{\left(  n\right)  }$ is the $n$-th derivate, and $g^{r}$ the $r$-th
exponential of $g$. For example we get for $k=0,1,2,\ldots$%
\begin{equation}
\left(  J^{-\frac{2k+1}{2}}\right)  ^{\left(  1\right)  }=-\frac{2k+1}%
{2}\epsilon J^{-\frac{2k+3}{2}} \label{Jr}%
\end{equation}
and integrating by parts we have%
\begin{equation}
\int\left(  J^{\frac{2k-1}{2}}\Psi\right)  =\frac{2\epsilon}{2k+1}\left(
J^{\frac{2k+1}{2}}\Psi-\int\left(  J^{\frac{2k-1}{2}}\Psi^{\left(  1\right)
}\right)  \right)  \label{intJr}%
\end{equation}
$\ $

It is suffice to prove that for any integer $k\geq0$ there exist $G_{k}\in
C^{\infty}\left(  I_{\varepsilon}\right)  $ and constant coefficients $a_{i}$
such that
\begin{equation}%
\begin{tabular}
[c]{c}%
$B_{k}=\sum_{i=0}^{k}\epsilon^{i+1}a_{i}\psi^{\left(  i\right)  }%
J^{-\frac{2\left(  k-i\right)  +1}{2}}$ and\\
$F^{\left(  k+1\right)  }=G_{k}+F_{k}$, where $F_{k}=B_{k}\int\left(
J^{\frac{2k-1}{2}}\psi^{\left(  k\right)  }\right)  $%
\end{tabular}
\ \label{Fk}%
\end{equation}

since $\lim_{t\rightarrow0}F_{k}\left(  t\right)  $ there exist
and it is finite. In fact by L'H\^{o}pital rule we see that
\[
m\geq n-1\geq0\Rightarrow\exists\lim_{t\rightarrow0}\frac{\int\left(
J^{m}\Psi\right)  }{J^{n}}=\lim_{t\rightarrow0}\frac{J^{m}\Psi}{nJ^{n-1}}%
\in\mathbb{R}%
\]
and (assuming (\ref{Fk})), this means that $\lim_{t\rightarrow0}F_{k}\left(
t\right)  \in\mathbb{R}$ when $k>0$.

Moreover for $k=0$
\begin{equation}
F^{\left(  1\right)  }=B_{0}\int\left(  J^{-\frac{1}{2}}\psi\right)
=F_{0}\text{, }(G_{0}=0,\text{ }B_{0}=2\epsilon J^{-\frac{1}{2}}\psi)
\label{F0}%
\end{equation}

Applying (\ref{intJr}) $\int\left(  J^{-\frac{1}{2}}\psi\right)
=2\epsilon\left(  J^{\frac{1}{2}}\psi-\int J^{\frac{1}{2}}\psi^{\left(
1\right)  }\right)  $ we see that
\[
F_{0}=4\psi^{2}-4J^{-\frac{1}{2}}\psi\int J^{\frac{1}{2}}\psi^{\left(
1\right)  }%
\]
by L'H\^{o}pital rule:%
\[
\lim_{t\rightarrow0}\frac{\int\left(  J^{\frac{1}{2}}\psi^{\left(  1\right)
}\right)  }{J^{\frac{1}{2}}}=\frac{1}{2}\lim_{t\rightarrow0}\frac{J^{\frac
{1}{2}}\psi^{\left(  1\right)  }}{J^{-\frac{1}{2}}}=\frac{1}{2}\lim
_{t\rightarrow0}J\psi^{\left(  1\right)  }=0
\]
and there exist $\lim_{t\rightarrow0}F_{0}\in\mathbb{R}$. Also this proves the
existence of $F_{k}$ and $G_{k}$ as in (\ref{Fk}) for $k=0$.

Assuming the existence of $B_{k}$ and $G_{k}$ as in (\ref{Fk}), we proceed by
induction. In order to construct $B_{k+1}$ we derive $F_{k}$ and we get:%
\[
F_{k}^{\left(  1\right)  }=B_{k}J^{\frac{2k-1}{2}}\psi^{\left(  k\right)
}+B_{k}^{\left(  1\right)  }\int\left(  J^{\frac{2k-1}{2}}\psi^{\left(
k\right)  }\right)
\]
a computation using (\ref{Jr}), (\ref{intJr}) gives for some constant
coefficients $b_{i}$ and $c_{i}$%
\[
B_{k}^{\left(  1\right)  }=-a_{0}\left(  \frac{2k+1}{2}\right)  J^{-\frac
{2k+3}{2}}\psi+\sum_{i=1}^{k+1}\epsilon^{i}b_{i}J^{\frac{2\left(  k-i\right)
+3}{2}}\psi^{\left(  i\right)  }%
\]%
\[
\int\left(  J^{\frac{2k-1}{2}}\psi^{\left(  k\right)  }\right)  =\frac
{2\epsilon}{2k+1}\left(  J^{\frac{2k+1}{2}}\psi^{\left(  k\right)  }%
-\int\left(  J^{\frac{2k+1}{2}}\psi^{\left(  k+1\right)  }\right)  \right)
\]%
\begin{align}
B_{k}^{\left(  1\right)  }\int\left(  J^{\frac{2k-1}{2}}\psi^{\left(
k\right)  }\right)   &  =-a_{0}\epsilon J^{-1}\psi\psi^{\left(  k\right)
}+\sum_{i=1}^{k+1}c_{i}\left(  \epsilon J\right)  ^{i-1}\psi^{\left(
i\right)  }\psi^{\left(  k\right)  }\label{sum1}\\
&  -\frac{2\epsilon}{2k+1}B_{k}^{\left(  1\right)  }\int\left(  J^{\frac
{2k+1}{2}}\psi^{\left(  k+1\right)  }\right) \nonumber
\end{align}%
\begin{equation}
B_{k}J^{\frac{2k-1}{2}}\psi^{\left(  k\right)  }=a_{0}\epsilon J^{-1}\psi
\psi^{\left(  k\right)  }+\sum_{i=1}^{k}a_{i}\left(  \epsilon J\right)
^{i-1}\psi^{\left(  i\right)  }\psi^{\left(  k\right)  } \label{sum2}%
\end{equation}
adding (\ref{sum1}) and (\ref{sum2}) we observe that cancel terms in $J^{-1}$
and we get for
\[
F_{k}^{\left(  1\right)  }=A_{k+1}+B_{k+1}\int J^{\frac{2k+1}{2}}\psi^{\left(
k+1\right)  }%
\]
where $A_{k+1}\in C^{\infty}\left(  I_{\varepsilon}\right)  $ and for some
constant coefficients $d_{i}$:%
\[
B_{k+1}=-\frac{2\epsilon}{2k+1}B_{k}^{\left(  1\right)  }=\sum_{i=0}%
^{k+1}\epsilon^{i+1}d_{i}J^{-\frac{2(k+1-i)+1}{2}}%
\]
an this end of the induction argument, and the proof.
\end{proof}

\begin{acknowledgement}
The author would like to thanks  \emph{Marek Kossowski} for explain me, by a
personal communication, the idea to establish the polar end metric definition,
which lives in the root of this paper, and my colleages \emph{Baldomero
Rubio}, since he is the author to the proof of the previous lemma, and
\emph{Eduardo Aguirre} by our conversations about the physical meaning of the
polar end spaces.
\end{acknowledgement}

\bibliographystyle{plain}

.
\end{document}